\documentclass[12pt]{article}
\usepackage{amsmath, amsthm, amsfonts, mathtools, bbm, amssymb, enumerate, hyperref, tikz}
\usetikzlibrary{plotmarks}

\newcommand{\G}{\ensuremath{\Gamma}}
\newcommand{\g}{\ensuremath{\gamma}}
\newcommand{\B}{\ensuremath{\beta}}
\newcommand{\e}{\ensuremath{\eta}}
\newcommand{\hg}{\ensuremath{\hat{\gamma}}}
\newcommand{\hb}{\ensuremath{\hat{\beta}}}
\newcommand{\N}{\ensuremath{\mathbb{N}}}
\newcommand{\K}{\ensuremath{\kappa}}
\newcommand{\hk}{\ensuremath{\hat{\kappa}}}

\newcommand{\A}{\ensuremath{\mathcal{A}}}

\newcommand{\z}{\ensuremath{\zeta}}
\newcommand{\hz}{\ensuremath{\hat{\zeta}}}
\newcommand{\R}{\ensuremath{\rho}}
\newcommand{\hr}{\ensuremath{\hat{\rho}}}
\renewcommand{\i}{{\rm i}}

\newtheorem{thm}{Theorem}[section]
\newtheorem{lem}[thm]{Lemma}

\newtheorem{rem}[thm]{Remark}

\everymath{\displaystyle}

\title{More on hypergeometric  L\'{e}vy processes}

\author{Emma L. Horton\thanks{Department of Mathematical Sciences, University of Bath, Claverton Down, Bath, BA2 7AY, UK. Email: \texttt{elh48@bath.ac.uk, a.kyprianou@bath.ac.uk}} $^{,}$\footnote{The idea for this work was proposed by Alexey Kuznetsov and it was conducted whilst EH was in receipt of a EPSRC funded SAMBa CDT undergraduate summer academic internship scholarship for which she is grateful.
}
\, and
Andreas E. Kyprianou\footnotemark[1]. 
}

\begin{document}

\maketitle

\begin{abstract}
Kuznetsov et al. \cite{Ketal} and  Kuznetsov and Pardo in~\cite{hg} 
 introduced the  family of Hypergeometric L\'evy processes. They appear naturally in the study of fluctuations of stable processes when one analyses stable processes through  the theory of positive self-similar Markov processes.  Hypergeometric L\'evy processes are defined through their  characteristic exponent, which, as a complex-valued function,  has four independent parameters.  Kyprianou et al. in~\cite{ehg} showed that the definition of a Hypergeometric L\'evy process could be taken to include a greater range of the aforesaid parameters than originally specified. In this short article, we push the parameter range even further. 
\end{abstract}

\section{Introduction}\label{c5_introduction}
Recall that a (killed) general one-dimensional L\'evy process is a stochastic process issued from the origin with stationary and independent increments and almost sure right continuous paths. We write $X=\{X_t : t\geq 0\}$ for its trajectory and $\mathbb{P}$ for its law. 
The law $\mathbb{P}$ of  a L\'evy process is characterized by its one-time transition probabilities. In particular there always exists a  quadruple  $(q, a, \sigma, \Pi)$ where $q\geq 0$ is the killing rate, $a\in\mathbb{R}$ is the linear coefficient, $\sigma\in \mathbb{R}$ is the Gaussian coefficient and $\Pi$ is a measure on $\mathbb{R}\backslash\{0\}$ satisfying $$\int_{\mathbb{R}}(1\wedge x^2)\Pi({\rm d}x)<\infty,$$ which gives the rate at which jumps of different sizes arrive,  such that 
\begin{equation}\label{eq1}
\mathbb{E}[{\rm e}^{\i z X_t}] = {\rm e}^{ t \psi(\i z)}, \qquad z \in\mathbb{R},
\end{equation}
where the Laplace exponent $\psi(z)$ is given by the L\'evy-Khintchine formula
\begin{equation}\label{def_psi}
\psi(z) = -q + a z  +\frac{1}{2}\sigma^2 z^2 + \int_{\mathbb{R}}\left( {\rm e}^{ z x} -1- z x \mathbf{1}_{(|x|<1)}\right)\Pi({\rm d}x).
\end{equation}

The {(spatial) Wiener--Hopf factorisation} of a L\'evy process $\xi$ with Laplace exponent $\psi$ consists of the equation
\[ 
\psi(z) = - \kappa(-z) \hat\kappa(z),\qquad  z \in \i \mathbb{R}, 
\]
where $\kappa$ and $\hat\kappa$ are the Laplace exponents of subordinators $H$ and $\hat H$, respectively, this time in the sense that $\mathbb{E}[e^{-\lambda H_1}\bigr] = e^{-\kappa(\lambda)}$ for $\Re \lambda \ge 0$. The subordinators $H$ and $\hat H$ are respectively known as the ascending and descending ladder height processes, and are related via a time-change to the running maximum and running minimum of the process $\xi$. The Wiener--Hopf factorisation has long been valued as the  insight into a great many fluctuation identities and, accordingly, has proved to be enormously important in the general theory of L\'evy processes. See the books \cite{Kyp} and \cite{bert} for more on the factorisation and the central role it plays in the analysis of a rich variety of problems.

There are relatively few instances where one finds concrete examples of L\'evy processes for which explicit identities exist for the characteristic exponent, the underlying quadruple $(q,a,\sigma, \Pi)$ as well as for the Wiener--Hopf factors. This is especially the case for L\'evy processes with two-sided jumps. One family of L\'evy processes for which this degree of tractability is available is the so-called {\it hypergeometric L\'evy} processes.  As we shall see below, hypergeometric L\'evy processes are  defined through their  characteristic exponent, which can be written in terms of four gamma functions and which accommodates for four independent parameters. These processes were introduced by Kuznetsov et al. \cite{Ketal} and  Kuznetsov and Pardo in~\cite{hg} and given a slightly more inclusive definition in \cite{ehg} by extending the parameter range. In this short article, we push the parameter range even further. 

\bigskip

We will now review the theory that has been developed for the hypergeometric L\'{e}vy processes in~\cite{hg} and~\cite{ehg} respectively. Note that the former reference refers to their generalisation of the L\'evy processes and in the latter reference as the {\it extended hypergeometric class}. As we are extending  the definition of this class even further, we prefer to  refer to all the L\'evy processes in this article as just hypergeometric L\'evy processes.

We say that  $\xi$ with parameters $(\B, \g, \hb, \hg)$ belongs to the hypergeometric class of L\'{e}vy processes if it has Laplace exponent which satisfies
\begin{equation}
\psi(z) = -\frac{\G(1-\B+\g-z)}{\G(1-\B-z)}\frac{\G(\hb+\hg+z)}{\G(\hb+z)}, \phantom{spaace} z \in \i\mathbb{R}.
\label{eq:1}
\end{equation}
In the existing literature hypergeometric class of L\'evy processes admits  parameter combinations in the domain $\mathcal{A}_1\cup\mathcal{A}_2$, where
\begin{equation*}
\A_1 = \{\B \le 1, \g \in (0,1), \hb \ge 0, \hg \in (0,1)\}
\end{equation*}
and 
\begin{equation*}
\A_2 = \{\B \in [1,2], \g, \hg \in (0,1), \hb \in [-1,0]; 1-\B+\hb+\g \ge 0, 1-\B+\hb+\hg \ge 0 \}.
\end{equation*}

The authors in~\cite{hg} also showed that, in the parameter regime $\mathcal{A}_1$,  the density of the L\'{e}vy measure is
\begin{equation}
\pi(x) = 
\begin{cases}
-\frac{\G(\e)}{\G(\e-\hg)\G(-\g)}{\rm e}^{-(1-\B+\g)x}{}_{2}{\mathcal{F}}_1(1+\g, \e; \e-\hg; {\rm e}^{-x}), &x>0, \\
-\frac{\G(\e)}{\G(\e-\g)\G(-\hg)}{\rm e}^{(\hb+\hg)x}{}_{2}{\mathcal{F}}_1(1+\hg, \e; \e-\g; {\rm e}^x), &x<0,
\end{cases} 
\label{eq:3}
\end{equation}
where ${}_{2}{\mathcal{F}}_1$ is the Hypergeometric function (which motivates the name of this class) and
\[
\e = 1-\B+\hb+\g+\hg.
\]
Moreover, it was shown that the Wiener--Hopf factors of a hypergeometric L\'evy process are
\begin{equation}
\K(\lambda) = \frac{\G(1-\B+\g+\lambda)}{\G(1-\B+\lambda)}, \phantom{spaace} \hk(\lambda) = \frac{\G(\hb+\hg+\lambda)}{\G(\hb+\lambda)},
\label{eq:2}
\end{equation}
for $\Re \lambda \geq 0$,
and that
in~\cite{ehg} it was shown that, in the parameter regime $\mathcal{A}_2$,  whilst the L\'evy density remains the same as in (\ref{eq:3}), the Wiener--Hopf factors are different. In that case, they can be written as
\begin{equation}
\K(\lambda) = (\lambda-\hb)\frac{\G(1-\B+\g+\lambda)}{\G(2-\B+\lambda)}, \phantom{spaace} \hk(\lambda) = (\B-1+\lambda)\frac{\G(\hb+\hg+\lambda)}{\G(1+\hb+\lambda)},
\label{eq:4}
\end{equation}
for $\Re \lambda\geq 0$.
%

\bigskip

One can derive the form of the L\'evy density in (\ref{eq:3}) in a relatively straightforward way for both parameter regimes $\mathcal{A}_1$ and $\mathcal{A}_2$,  by using  the theory of meromorphic L\'{e}vy processes (see~\cite{mero}). Recall that a L\'evy process belongs to the meromorphic class if its L\'evy measure is absolutely continuous with respect to Lebesgue measure and the density takes the form 
\begin{equation}
\pi(x) = \mathbf{1}_{\{x>0\}} \sum\limits_{n\ge 1} a_n \rho_n e^{-\rho_n x}+ \mathbf{1}_{\{x<0\}} \sum\limits_{n\ge 1} \hat a_n \hat \rho_n e^{\hat \rho_n x}, \qquad x\in \mathbb{R},
\label{mmphic}
\end{equation}
 where all the coefficients  $a_n$, $\hat a_n$, $\rho_n$, $\hat \rho_n$ are positive, 
the sequences $\{\rho_n\}_{n \ge 1}$ and $\{\hat \rho_n\}_{n\ge 1}$ are strictly increasing, and $\rho_n \to +\infty$ and $\hat \rho_n \to +\infty$ as $n\to +\infty$.
Part of the  theory of meromorphic L\'evy processes relies on the fact that $\{\z_k, -\hz_k \colon k \in \N \}$ and $\{\R_k, -\hr_k \colon k \in \N \}$ are, respectively, the roots and poles of the Laplace exponent $\psi$, which necessarily interlace in the following sense
\begin{equation}
\dots < -\hr_2 < - \hz_2 < -\hr_1 < -\hz_1 < 0 < \z_1 < \R_1 < \z_2 < \R_2 < \dots
\label{eq:5}
\end{equation}

For the case of hypergeometric L\'evy processes, since $\G(z)$ has simple poles at the points $\{-k \colon k \in \N_0\}$, the exponent $\psi$ in (\ref{eq:1}) has simple poles at the points $\{1-\B+\g+k, -\hb-\hg-k \colon k \in \N_0\}$ and simple roots at the points $\{1-\B+k, -\hb-k \colon k \in \N_0\}$.
As an example, consider a hypergeometric process such that  $\B=0.9$, $\g=0.5$, $\hb=0.2$ and $\hg=0.3$. Note that these parameters lie in $\mathcal{A}_1$. Then $\psi$ has simple poles of the form $\{-0.5-k, 0.6+k \colon k \in \N_0\}$ and simple roots of the form $\{-0.2-k, 0.1+k \colon k \in \N_0\}$. The roots (dots) and poles (crosses) interlace as shown in the graph below.
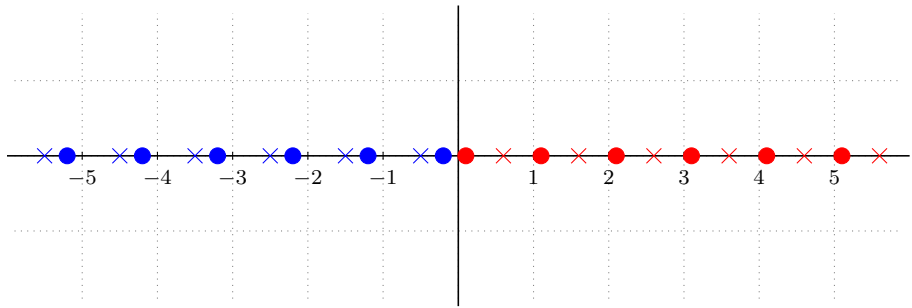
\begin{center}
\begin{tikzpicture}
\draw[step=1cm,gray,dotted] (-5.9,-1.9) grid (5.9,1.9);
\draw[semithick] (-6,0) -- (6,0);
\draw[semithick] (0,-2) -- (0,2);

\foreach \x in {-5,-4,-3,-2,-1,1,2,3,4,5}
	\draw (\x cm, 1.5pt) -- (\x cm, -1.5pt) node[anchor=north] {\scriptsize $\x$};

\draw [blue] plot [only marks, mark=x, mark size=4] coordinates {(-5.5,0) (-4.5,0) (-3.5,0) (-2.5,0) (-1.5,0) (-0.5,0)}; 
\draw [blue] plot [only marks, mark=*, mark size=3] coordinates {(-5.2,0) (-4.2,0) (-3.2,0) (-2.2,0) (-1.2,0) (-0.2,0)}; 
\draw [red] plot [only marks, mark=x, mark size=4] coordinates {(0.6,0) (1.6,0) (2.6,0) (3.6,0) (4.6,0) (5.6,0)}; 
\draw [red] plot [only marks, mark=*, mark size=3] coordinates {(0.1,0) (1.1,0) (2.1,0) (3.1,0) (4.1,0) (5.1,0)}; 
\end{tikzpicture}
\end{center}
The red roots and poles emerge from the factor 
\begin{equation}
\frac{\G(1-\B+\g-z)}{\G(1-\B-z)}
\label{red}
\end{equation}
 and the blue roots and poles emerge from the factor 
 \begin{equation}
 \frac{\G(\hb+\hg+z)}{\G(\hb+z)}.
 \label{blue}
 \end{equation}
For the parameter choices here, the ascending and descending Wiener--Hopf factors are determined entirely by the red and blue roots and poles. Said another way, the two factors are precisely (\ref{red}) and (\ref{blue}) respectively.

\bigskip

Now consider another hypergeometric process whose parameters are  $\B=1.1$, $\g=0.5$, $\hb=-0.2$ and $\hg=0.8$, i.e. belonging to the regime $\mathcal{A}_2$. We now have simple poles of $\psi$ in (\ref{eq:1}) of the form $\{-0.6-k, 0.4+k \colon k \in \N_0\}$ and simple roots of the form $\{0.2-k, -0.1+k \colon k \in \N_0\}$. Again, the interlacement of these roots (dots) and poles (crosses) are represented graphically below.  Once again, the red roots and poles emerge from the factor (\ref{red}) and the blue roots and poles from the factor (\ref{blue}). 
Note, however, for these parameter choices,  the ascending and descending Wiener--Hopf factors are not determined by (\ref{red}) and (\ref{blue}) respectively. 
It is the positive roots and poles which determine the ascending factor and the negative roots and poles which determine the descending factor. In this sense, the two Wiener--Hopf factors are given by (\ref{eq:4}).

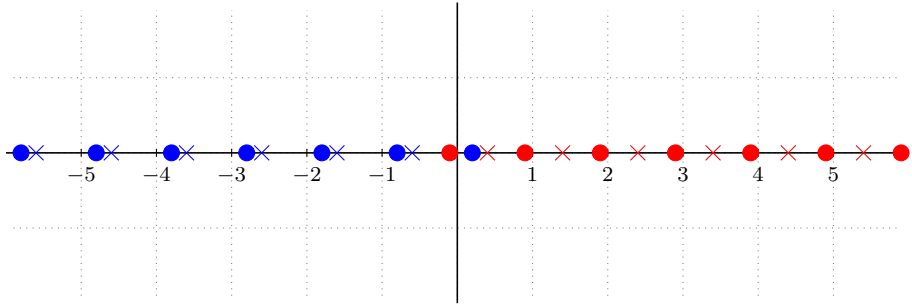
\begin{center}
\begin{tikzpicture}
\draw[step=1cm,gray,dotted] (-5.9,-1.9) grid (5.9,1.9);
\draw[semithick] (-6,0) -- (6,0);
\draw[semithick] (0,-2) -- (0,2);

\foreach \x in {-5,-4,-3,-2,-1,1,2,3,4,5}
	\draw (\x cm, 1.5pt) -- (\x cm, -1.5pt) node[anchor=north] {\scriptsize $\x$};
		
\draw [blue] plot [only marks, mark=x, mark size=4] coordinates {(-5.6,0) (-4.6,0) (-3.6,0) (-2.6,0) (-1.6,0) (-0.6,0)}; 
\draw [blue] plot [only marks, mark=*, mark size=3] coordinates {(-5.8,0) (-4.8,0) (-3.8,0) (-2.8,0) (-1.8,0) (-0.8,0) (0.2,0)}; 
\draw [red] plot [only marks, mark=x, mark size=4] coordinates {(0.4,0) (1.4,0) (2.4,0) (3.4,0) (4.4,0) (5.4,0)}; 
\draw [red] plot [only marks, mark=*, mark size=3] coordinates {(-0.1,0) (0.9,0) (1.9,0) (2.9,0) (3.9,0) (4.9,0) (5.9,0)}; 
\end{tikzpicture}
\end{center}
\section{Extending the parameter domain further}\label{c5_eehg}

We first define the two sets of admissible parameters for that we would like to include in the definition of hypergeometric L\'evy processes:
\begin{equation*}
\A_3 = \bigcup_{n=0}^\infty \{\B \in [0,1], \g, \hg \in (0,1), \hb \in [-(n+1), -n]; 1-\B+\hb+\hg+n \le 0, 1-\B+\hb+\g+n \ge 0 \}
\end{equation*}
and
\begin{equation*}
\A_4 = \bigcup_{n=1}^\infty \{\hb \in [0,1], \g,\hg \in (0,1), \B \in [n,n+1]; n-\B+\hb+\hg \ge 0, n-\B+\hb+\g \le 0 \}
\end{equation*}
Our objective is to show that the the expression in (\ref{eq:1}) is a characteristic exponent of  a L\'evy process when $(\beta, \gamma, \hat\beta, \hat\gamma)\in \mathcal{A}_3\cup\mathcal{A}_4$.
We shall preemptively consider two examples from these two parameter regimes and see how this suggests we might be able to extend the parameter range of hypergeometric L\'evy processes.

For the first example, the parameters take values in $\A_3$. Fix $\B=0.8$, $\g=0.4$, $\hb=-1.9$ and $\hg=0.1$.
We have simple roots of $\psi$ in (\ref{eq:1}) at $\{1.5-k, 0.2+k \colon k \in \N_0\}$ and simple poles at $\{1.4-k, 0.6+k \colon k \in \N_0\}$. Notice that we are in the case when $n$ (in the definition of $\mathcal{A}_3$) is equal to 1. The interlacement of the roots and poles still holds, see the figure below. Once again, we have coloured red the roots and poles that come from (\ref{red}) and blue the roots and poles that come from (\ref{blue}).
 By allowing $\hb$ to be pushed further into the negative reals, we see root-pole pairs from (\ref{blue}) have passed over to the positive side of the origin, but all the root-pole pairs from (\ref{red}) remain positive. Clearly (\ref{red}) and (\ref{blue}) cannot represent the ascending and descending Wiener--Hopf factors respectively. However, in the spirit of (\ref{eq:4}), by trading linear factors between the two expressions in (\ref{red}) and (\ref{blue}) we can associate all of the positive and negative root-poles pairs to two different factors. As we shall see, this turns out to identify a Wiener--Hopf factorisation of the form
 \[
 \psi(\i z) = \left[\frac{\G(1-\B+\g-\i z)}{\G(1-\B-\i z)}\frac{(-\hb-\i z)}{(-\hb-\hg-\i z)}\frac{(-\hb-1-\i z)}{(-\hb-\hg-1-\i z)}
\right]\times \left[\frac{\G(2+\hb+\hg+\i z)}{\G(2+\hb+\i z)}, 
\right]
 \]
 for $z\in\mathbb{R}$.
  
\begin{center}
\begin{tikzpicture}[scale=1.5]
\draw[step=1cm,gray,dotted] (-3.9,-0.9) grid (4.9,0.9);
\draw[semithick] (-4,0) -- (5,0);
\draw[semithick] (0,-1) -- (0,1);

\foreach \x in {-3,-2,-1,1,2,3,4}
	\draw (\x cm, 1.5pt) -- (\x cm, -1.5pt) node[anchor=north] {\scriptsize $\x$};
		
\draw [blue] plot [only marks, mark=x, mark size=2.5] coordinates {(-3.6,0) (-2.6,0) (-1.6,0) (-0.6,0) (0.4,0) (1.4,0)}; 
\draw [blue] plot [only marks, mark=*, mark size=1.5] coordinates {(-3.5,0) (-2.5,0) (-1.5,0) (-0.5,0) (0.5,0) (1.5,0)}; 
\draw [red] plot [only marks, mark=x, mark size=2.5] coordinates {(0.6,0) (1.6,0) (2.6,0) (3.6,0) (4.6,0)}; 
\draw [red] plot [only marks, mark=*, mark size=1.5] coordinates {(0.2,0) (1.2,0) (2.2,0) (3.2,0) (4.2,0)}; 
\end{tikzpicture}
\end{center}

For the second example, let us consider the case that  our parameters take values in $\A_4$. Fix $\B=2.7$, $\g=0.2$, $\hb=0.3$ and $\hg=0.5$. We have now pushed $\B$ further into the positive reals, which corresponds to passing root-pole pairs over to the negative side of 0 as seen below. In a similar way to the previous, we can trade linear factors between (\ref{red}) and (\ref{blue}) to ensure that all the positive and negative root-pole pairs appear in different factors of $\psi$ in (\ref{eq:1}). The obvious candidate for the Wiener--Hopf factorisation would then take the form 
\[
\psi(\i z) = \left[ \frac{\G(3-\B+\g-\i z)}{\G(3-\B-\i z)}
\right]\times\left[\frac{\G(\hb+\hg+\i z)}{\G(\hb+\i z)} \frac{(\B-1+\i z)}{(\B-\g-1+\i z)}\frac{(\B-2+\i z)}{(\B-\g-2+\i z)}\right]
\]

\begin{center}
\begin{tikzpicture}[scale=1.5]
\draw[step=1cm,gray,dotted] (-5.9,-0.9) grid (3.9,0.9);
\draw[semithick] (-6,0) -- (4,0);
\draw[semithick] (0,-1) -- (0,1);

\foreach \x in {-5,-4,-3,-2,-1,1,2,3}
	\draw (\x cm, 1.5pt) -- (\x cm, -1.5pt) node[anchor=north] {\scriptsize $\x$};
		
\draw [blue] plot [only marks, mark=x, mark size=2.5] coordinates {(-5.8,0) (-4.8,0) (-3.8,0) (-2.8,0) (-1.8,0) (-0.8,0)}; 
\draw [blue] plot [only marks, mark=*, mark size=1.5] coordinates {(-5.3,0) (-4.3,0) (-3.3,0) (-2.3,0) (-1.3,0) (-0.3,0)}; 
\draw [red] plot [only marks, mark=x, mark size=2.5] coordinates {(-1.5,0) (-0.5,0) (0.5,0) (1.5,0) (2.5,0) (3.5,0)}; 
\draw [red] plot [only marks, mark=*, mark size=1.5] coordinates {(-1.7,0) (-0.7,0) (0.3,0) (1.3,0) (2.3,0) (3.3,0)}; 
\end{tikzpicture}
\end{center}

\begin{thm}\label{c5_main_result_1}
There exists a L\'{e}vy process $\xi$ with Laplace exponent $\psi$ as given in~\eqref{eq:1}, where $(\B, \g, \hb, \hg)$ take values in
\begin{enumerate}
\item[(i)] \label{1} $\A_1$. Its Wiener--Hopf factors are given by~\eqref{eq:2} and its L\'{e}vy density is given by~\eqref{eq:3}.
\item[(ii)] \label{2} $\A_2$. Its Wiener--Hopf factors are given by~\eqref{eq:4} and its L\'{e}vy density is given by~\eqref{eq:3}.
\item[(iii)] \label{3} $\A_3$. Its Wiener--Hopf factors are given by 
\begin{equation}
\K(\lambda) = \frac{\G(1-\B+\g+\lambda)}{\G(1-\B+\lambda)}\prod_{j=0}^n\frac{(-\hb-j+\lambda)}{(-\hb-\hg-j+\lambda)}, \qquad \lambda \geq 0.\label{plane}
\end{equation}
and 
\begin{equation*}
\hk(\lambda) = \frac{\G(n+1+\hb+\hg+\lambda)}{\G(n+1+\hb+\lambda)}, \qquad \lambda \geq 0.
\end{equation*}
Its L\'{e}vy density is given by
\begin{equation*}
\pi(x) =
\begin{cases}
-\frac{\G(\e)}{\G(-\hg)\G(\e-\g)}{\rm e}^{(\hb+\hg)x}\sum_{k=0}^n\frac{(1+\hg)_k(\e)_k}{(\e-\g)_k}\frac{{\rm e}^{kx}}{k!}\\
\phantom{spaace} -\frac{\G(\e)}{\G(\e-\hg)\G(-\g)}{\rm e}^{-(1-\B+\g)x} {}_{2}{\mathcal{F}}_1(1+\g, \e; \e -\hg; {\rm e}^{-x}), &x > 0, \\
-\frac{\G(\e)}{\G(-\hg)\G(\e -\g)}{\rm e}^{(\hb+\hg)x}\sum_{k \ge n+1} \frac{(1+\hg)_k (\e)_k}{(\e -\g)_k}\frac{{\rm e}^{kx}}{k!}, &x<0.
\end{cases}
\end{equation*}
\item[(iv)] \label{4} $\A_4$. Its Wiener--Hopf factors are given by 
\begin{equation*}
\K(\lambda) = \frac{\G(1+n-\B+\g+\lambda)}{\G(1+n-\B+\lambda)}, \qquad \lambda\geq 0
\end{equation*}
and
\begin{equation*}
\hk(\lambda) = \frac{\G(\hb+\hg+\lambda)}{\G(\hb+\lambda)} \prod_{j=1}^n \frac{\B-j+\lambda}{\B-\g-j+\lambda}, \qquad \lambda\geq 0.
\end{equation*}
Its L\'{e}vy density is given by
\begin{equation}
\pi(x) =
\begin{cases}
-\frac{\G(\e)}{\G(-\g)\G(\e-\hg)}{\rm e}^{-(1-\B+\g)x}\sum_{k \ge n}\frac{(1+\g)_k (\e)_k}{(\e-\hg)_k}\frac{{\rm e}^{-kx}}{k!}, &x>0 \\
-\frac{\G(\e)}{\G(-\g)\G(\e-\hg)}{\rm e}^{-(1-\B+\g)x}\sum_{k=0}^{n-1}\frac{(1+\g)_k(\e)_k}{(\e-\hg)_k}\frac{{\rm e}^{-kx}}{k!} \\
\phantom{spaaaace} -\frac{\G(\e)}{\G(-\hg)\G(\e-\g)}{\rm e}^{(\hb+\hg)x}{}_{2}{\mathcal{F}}_1(1+\hg, \e;\e-\g;{\rm e}^x), &x<0
\end{cases}
\label{A4density}
\end{equation}
\end{enumerate}
For all parameter sets, the process is killed at rate 
\begin{equation*}
q=\frac{\G(1-\B+\g)}{\G(1-\B)}\frac{\G(\hb+\hg)}{\G(\hb)}
\end{equation*}
Furthermore, the process $\xi$ has no Gaussian component. When $\g+\hg \in (0,1)$ \textup(resp. $\g+\hg \in [1,2)$\textup) $\xi$ has paths of bounded variation and no drift \textup(resp. paths of unbounded variation\textup).
\end{thm}
\begin{rem}\rm
Note that the parameter ranges $\A^{\hb}_{EHG}$ and $\A^{\B}_{EHG}$ defined in~\cite[Remark 3]{ehg} are included in the parameter ranges $\A_3$ and $\A_4$ respectively.
\end{rem}
\begin{rem}\rm\label{prop1}
Let $\xi$ be a L\'{e}vy process as defined in~\eqref{c5_main_result_1} whose parameters lie in $\A_4$. By noting that the ascending Wiener--Hopf factor belongs to the class of $\beta$-subordinators, we have from standard literature (cf. \cite{Ketal}) and direct computation, with the help of the beta integral, that  the L\'{e}vy density of the ascending ladder height $\K$ is
\begin{equation*}
\nu(x) = \frac{\g}{\G(1-\g)}(1-{\rm e}^{-x})^{-\g-1} {\rm e}^{(n+1-\B+\g)x}, \phantom{spaace} x>0.
\end{equation*}
The density for its potential measure is given by
\begin{equation*}
{u}(x) = \frac{1}{\G(\g)} {\rm e}^{-(n+1-\B)x} (1-{\rm e}^{-x})^{\g-1}, \phantom{spaace} x \ge 0.
\end{equation*}
The same objects for the descending ladder exponent $\hat\kappa$ are somewhat more complicated. Its potential cannot be written so easily in closed form, however, in principle, since $\psi$ is a meromorphic L\'evy process, it can be written as a mixture of exponentials. 
\end{rem}

\section{Proof of Theorem~\ref{c5_main_result_1}}\label{c5_proofs}
For the proof of (i), we refer the reader to~\cite{hg}. Also, the proof of (ii) can be found in~\cite{ehg}. We will now prove (iii). Our strategy will be to first show that the claimed factors are indeed Bernstein functions (i.e. Laplace exponents of subordinators) and that they belong together in a Wiener--Hopf factorisation, i.e. the function $\psi$ in (\ref{eq:1}) is a characteristic exponent of a L\'evy process.

We will first work under the assumption that $1-\B+\hb+\hg+n < 0$ and $1-\B+\hb+\g+n > 0$. Recall that the function $\psi$ has simple poles at the points $\{1-\B+\g+k, -\hb-\hg-k \colon k \in \N_0\}$ and simple roots at the points $\{1-\B+k, -\hb-k \colon k \in \N_0\}$.

For a given choice of parameters $(\B, \g, \hb, \hg) \in \A_3$, there are $n+1$ positive poles, say $\R_j^*$, of the form $-\hb-\hg-j$ for $j=0,\dots,n$ with the rest of the positive poles, say $\R_l^{**}$, being of the form $-\B+\g+l$ for $l \in \N$. Similarly, there are $n+1$ positive roots, say $\z_j^*$, of the form $-\hb-j$ for $j = 0,\dots,n$ with the rest of the positive roots, say $\z_l^{**}$, being of the form $-\B+l$ for $l \in \N$. The conditions $1-\B+\hb+\hg < 0, 1-\B+\hb+\g >0$, which are included in $\A_3$, mean that the positive poles and roots interlace as follows: 
\begin{equation*}
0 < \z_1^{**} < \R_0^* < \z_0^* < \R_1^{**} < \z_2^{**} < \R_1^* < \dots < \z_n^* < \R_{n+1}^{**} < \z_{n+2}^{**} < \R_{n+2}^{**} < \dots
\end{equation*}
For $k \in \N$, denote the ordered (i.e. in order of interlacement) positive poles and roots by $\R_k$ and $\z_k$ respectively, then, for $z\in\mathbb{C}$ such that the left hand side is well defined,
\begin{align}
\prod_{k \ge 1} \frac{1+z/\z_k}{1+z/\R_k} &= \prod_{l \ge 1}\frac{1+z/\z_l^{**}}{1+z/\R_l^{**}}\prod_{j=0}^n\frac{1+z/\z_j^*}{1+z/\R_j^*} \notag\\
&= \prod_{l \ge 1}\frac{1+\frac{z}{-\B+l}}{1+\frac{z}{-\B+\g+l}}\prod_{j=0}^n\frac{1+\frac{z}{-\hb-j}}{1+\frac{z}{-\hb-\hg-j}} \label{derivedlater} \\
&\approx \frac{\G(1-\B+\g+z)}{\G(1-\B+z)}\prod_{j=0}^n\frac{-\hb-j+z}{-\hb-\hg-j+z}\notag\\
&=\vcentcolon \K(z)\notag,
\end{align}
where $\approx$ means ``up to a constant''.
Next, we need to apply Lemma 1 of \cite{mero} (see also equations (16)  and (17) in that paper). For convenience, we reproduce it below in paraphrased form.

\begin{lem}\label{import}Assume that we have two increasing sequences $\rho=\{\rho_n\}_{n\ge 1}$ and $\zeta=\{\zeta_n\}_{n \ge 1}$ of positive numbers, such that
 $\rho_n\to +\infty$ as $n\to +\infty$ and the following interlacing condition is satisfied:
 \[
 \zeta_1< \rho_1 < \zeta_2 < \rho_2 < ...
\]
 Define
 \[
\phi(z)= \prod\limits_{n\ge 1}  \frac{1+{z}/{\rho_n}}{1+{z}/{\zeta_n}}, \;\;\; z>0.
\]
Then, for all $z>0$,
\begin{eqnarray*}
 \phi(z)&=& \texttt{a}_0(\rho,\zeta) + \int_{\mathbb{R}^+} \left[ \sum\limits_{n\ge 1} \texttt{a}_n(\rho,\zeta) \zeta_n e^{-\zeta_n x} \right] e^{-zx} {\rm d} x,\\
\frac{1}{ \phi(z)}&=& 1+z \texttt{b}_0(\zeta,\rho) +\int_{\mathbb{R}^+} \left[\sum\limits_{n\ge 1} \texttt{b}_n(\zeta,\rho) 
 \rho_n e^{-\rho_n x} \right] \left(1-e^{-zx}  \right) {\rm d} x,
\label{see-the-drift}
\end{eqnarray*}
where 
 \begin{eqnarray*}
 \texttt{a}_0(\rho,\zeta)&=&\lim_{n\to +\infty}\prod\limits_{k=1}^n \frac{\zeta_k}{\rho_k}, \;\;\;
\texttt{a}_n(\rho,\zeta)=\left(1-{\zeta_n}/{\rho_n}\right) \prod\limits_{\substack{k\ge 1 \\ k\ne n}}  \frac{1-{\zeta_n}/{\rho_k}}{1-{\zeta_n}/{\zeta_k}}, \\
 \texttt{b}_0(\zeta,\rho)&=&\frac{1}{\zeta_1}\lim_{n\to +\infty}\prod\limits_{k=1}^n \frac{\rho_k}{\zeta_{k+1}}, \;\;\;
\texttt{b}_n(\zeta,\rho)=-\left(1-\frac{\rho_n}{\zeta_n}\right) \prod\limits_{\substack{k\ge 1 \\ k\ne n}} 
 \frac{1-{\rho_n}/{\zeta_k}}{1-{\rho_n}/{\rho_k}}.
 \end{eqnarray*}
 Moreover, $\texttt{a}_0(\rho,\zeta)\ge 0$, $\texttt{b}_0(\zeta,\rho)\ge 0$ and for all $n\ge 1$ we have $\texttt{a}_n(\rho,\zeta) > 0$, $\texttt{b}_n(\zeta,\rho) > 0$.

\end{lem}

We see that, for $\lambda\geq 0$, $\kappa(\lambda)$ is a Bernstein function (in the notation of Lemma  \ref{import}, $\phi(z) =1/\kappa(z)$) and so is $\lambda/\kappa(\lambda)$, $\lambda\geq 0$. Recall that $\kappa(\lambda)$, $\lambda\geq 0$  is a special Bernstein function if and only if $\lambda /\kappa(\lambda)$, $\lambda \geq 0$ is a Bernstein function, in which case the underlying L\'evy measure associated to either of these two is absolutely continuous with non-increasing density. We thus have that  $ \kappa(\lambda)$ is a special Bernstein function.

Now denote the negative poles by $\hr_k$. For $k \ge 1$ they take the form $\hb+\hg+k+n$. Similarly, the zeros $\hz_k$ take the form $\hb+k+n$. Note that these poles and roots interlace in the following way
\begin{equation*}
\dots < -\hr_2 < -\hz_2 < -\hr_1 < -\hz_1 < 0
\end{equation*}
Therefore
\begin{align*}
\prod_{k \ge 1} \frac{1+z/\hz_k}{1+z/\hr_k} &= \prod_{k \ge 1}\frac{1 + \frac{z}{n+\hb+\hg+k}}{1+\frac{z}{n+\hb+k}} \\
&\approx \frac{\G(n+1+\hb+\hg+z)}{\G(n+1+\hb+z)} \\
&=\vcentcolon \hk(z)
\end{align*}
A similar argument to the one in the previous paragraph shows similarly that $\hat{\kappa}(\lambda)$ is also a special Bernstein function. 
%

Next, using the relation $\G(x+1)=x\G(x)$ repeatedly, one easily verifies  that, for the expression given in (\ref{eq:1}), $\psi(z) = -\kappa(-z)\hat{\kappa}(z)$. We can now appeal to Vigon's theory of philanthropy \cite[Chapter 7]{Vigon} (see also Section 6.6 of \cite{Kyp}), which shows that $\psi(z)$ is the Laplace exponent of a L\'{e}vy process. Specifically Vigon's theory of philanthropy states that if two subordinators have Laplace exponents $\phi_i(\lambda)$, $\lambda \geq 0$, $i=1,2$, each of which has absolutely continuous and non-increasing L\'evy density, then they belong together in a Wiener--Hopf factorisation (they become `friends'). That is to say, $\phi_1(-\i \theta)\phi_2(\i \theta)$, $\theta\in\mathbb{R}$ is the characteristic exponent of a L\'evy process.

To show $\xi$ is a meromorphic L\'{e}vy process, apply \cite[Theorem 1(v)]{mero} in the killed case and \cite[Corollary 2]{mero} in the unkilled case. For L\'{e}vy processes in the meromorphic class, it is known that their L\'{e}vy measure has a density, $\pi$, of the form given in (\ref{mmphic}).
We will now compute the coefficients $a_k\R_k$ and $\hat{a}_k\hr_k$ in the representation (\ref{mmphic}) and accordingly prove that $\pi$ is equivalent to the expression given in the statement of the theorem.\\
Firstly, label the $\R_k$ as follows: 
\begin{alignat*}{2}
\R_k &= -\hb-\hg-k+1, &\quad 1 \le &k \le n+1 \\
\R_k &= -\B+\g+k-n-1, & &k \ge n+2
\end{alignat*}
i.e. the first $n+1$ terms correspond to the poles that have been derived from (\ref{blue}) and the second group of terms correspond to the poles that have been derived from (\ref{red}) in the regime $\A_1$.\\
\noindent
Using the relation 
\begin{equation*}
a_k \R_k = -{\rm Res}(\psi(z) \colon z = \R_k)
\end{equation*}
we find that for $1 \le k \le n+1$
\begin{equation*}
a_k \R_k = -\frac{(-1)^{k-1}}{(k-1)!}\frac{1}{\G(1-\hg-k)}\frac{\G(\e -1+k)}{\G(\e -\g-1+k)}
\end{equation*}
and for $k \ge n+2$
\begin{equation*}
a_k \R_k = -\frac{(-1)^{k-n-2}}{(k-n-2)!}\frac{1}{\G(2-\g -k+n)}\frac{\G(\e -2-n+k)}{\G(\e -2-\hg -n+k)}
\end{equation*}
Now consider the negative poles. For $k \ge 1$ we can write
\begin{equation*}
\hr_k = \hb+\hg+k+n
\end{equation*}
and using 
\begin{equation*}
\hat{a}_k\hr_k = {\rm Res}(\psi(z) \colon z = -\hr_k)
\end{equation*}
we can compute the coefficients $\hat{a}_k\hr_k$ in a similar fashion.
Substituting these expressions into (\ref{mmphic}) and recalling that 
\[
{}_{2}{\mathcal{F}}_1(1+\hg, \e;\e-\g;z) = \sum_{n = 0}^\infty \frac{(1+\hg)_n(\e)_n}{(\e-\g)_n}\frac{ z^n}{n!}, \qquad z\geq 0,
\]
we obtain the density of the L\'{e}vy measure as stated in the theorem. 

\bigskip

Now suppose at least one of the following two equalities holds 
\begin{equation}
1-\B+\hb+\hg+n = 0, \quad 1-\B+\hb+\g+n = 0.
\label{eq:7}
\end{equation}
Then a finite number of the roots and poles are equal so cancel each other out in the product representation of the Wiener-Hopf factors. The remaining roots and poles still interlace in the required fashion so all of the expressions previously calculated still hold. However, this means that the coefficients corresponding to the poles that have vanished are no longer present in the calculation for $\pi$.\\ 

One can prove part (iv) of the proposition in a similar way to part (iii) and we omit this proof for the sake of brevity.
The expression for the killing rate follows by noting that $q=-\psi(0)$.

The proof of the claim about bounded and unbounded variation follows along the same lines as the proof given in~\cite[Proposition 1]{ehg}, however, we will say a few words about it. Formula 8.328.1 in~\cite{tables} implies that
\begin{equation}
\psi(iz) = O(\lvert z \rvert^{\g+\hg}), \phantom{spaace} z \to \infty, z \in \mathbb{R}.
\label{eq:8}
\end{equation}
Proposition 2 from~\cite{bert} then implies that $\xi$ has no Gaussian component. For the claims about the bounded and unbounded variation, the same proof as in~\cite[Proposition 1]{ehg} still applies for all parameter regimes. For example, if we consider a process in $\A_4$, the finite sum given in the L\'{e}vy density for $x<0$ is of constant order, and since the density for $x>0$ is a hypergeometric function minus this finite sum, the same asymptotic behaviour still holds for $\pi(x)$. Finally, when $\xi$ has paths of bounded variation,~\cite[Proposition 2]{bert} and~\eqref{eq:8} shows that $\xi$ has no drift.
\hfill$\square$

\end{document}